# Asymptotics of sums of functions of primes located on an arithmetic progression

Victor Volfson

ABSTRACT In the paper, we investigate the problem of the distribution of sums of functions of prime numbers located on an arithmetic progression. This problem is closely related to the problem of the distribution of prime numbers on an arithmetic progression. Based on this distribution, a general formula was obtained for the asymptotic estimate of the sums of functions of primes, and also asymptotics was found for the sums of various functions of primes on a geometric progression. Several assertions about asymptotics of sums of functions of prime numbers on a geometric progression are proved. Necessary and sufficient conditions for the existence of these asymptotics are also proved.

Keywords: asymptotic estimate for the number of primes on a geometric progression, generalized Riemann conjecture, sum of functions of primes, Abel's summation formula, asymptotic estimate for the remainder, necessary and sufficient conditions for the existence of the asymptotics.

## 1. INTRODUCTION

We investigated the problem of the distribution of sums of functions of prime numbers in the natural series in [1], which is directly related with the problem of the distribution of primes in the natural series.

We will consider the problem of the distribution of the sums of functions of prime numbers that are on an arithmetic progression in the paper. Naturally, this problem is related to the problem of the distribution of prime numbers on an arithmetic progression.

It is proved [2] that for all $k, l$ in geometric progression $k, l$, where $(k, l) = 1$, the asymptotic of the number of primes on the ray $[2, x)$ is:

$$\pi_l(k, x) = \frac{1}{\varphi(k)} \int_2^x \frac{dt}{\log t} + O(\frac{x}{e^{c\sqrt{\log(x)}}}), \qquad (1.1)$$

where $c > 0$ is constant.

Based on (1.1), one can obtain a coarser asymptotic estimate for the number of primes located on a geometric progression $kn + l$ on a ray $[2, x)$.

Integrating by parts, we get:

$$\int_2^x \frac{du}{\log(u)} = \frac{u}{\log(u)}\Big|_2^x - \int_2^x u\, d(\frac{1}{\log(u)}) = \frac{x}{\log(x)} - \frac{2}{\log 2} + \int_2^x \frac{du}{\log^2 u} = \frac{x}{\log(x)} + O(\frac{x}{\log^2 x})$$

Substituting the resulting expression in (1.1), we get the asymptotic:

$$\pi_l(k, x) = \frac{x}{\varphi(k)\log(x)} + O(\frac{x}{\log^2 x}) + O(\frac{x}{e^{c\sqrt{\log(x)}}}) = \frac{x}{\varphi(k)\log(x)} + O(\frac{x}{\log^2 x}), \qquad (1.2)$$

which sometimes use.

Expression (1.1) of the asymptotic law of distribution of primes was first obtained in 1896 by Hadamard and Valais-Poussin [2] based on an estimate of the boundary of the region of absence of nontrivial zeros of the Riemann zeta function:

$$\sigma \geq 1 - \frac{c}{\log T}, \qquad (1.3)$$

where $T \geq 10$.

Vinogradov [3] proved a new bound for the boundary of the region of absence of nontrivial zeros of the Riemann zeta function in 1958:

$$\sigma \geq 1 - \frac{c}{\log T^{2/3-\xi}}. \qquad (1.4)$$

Having in mind (1.4), a new asymptotic estimate for the remainder in formula (1.1) was obtained:

$$O\left(\frac{x}{e^{c(\log(x))^{0,6-\xi}}}\right). \qquad (1.5)$$

The result (1.5) is practically not improved to this day. Meanwhile, based on Riemann conjecture [4], which has not yet been proved, all non-trivial zeros of the zeta function are on the straight line $\sigma = 1/2$.

Like the Riemann conjecture, which has a corollary on the distribution of primes in natural numbers, the generalized Riemann conjecture has a corollary on the distribution of primes in an arithmetic progression.

Dirichlet character is a fully multiplicative arithmetic function χ such that there is a positive integer k with χ (n + k) = χ (n) for all n and χ (n) = 0 if (n, k)> 1. If there is such a character, then we define the corresponding Dirichlet L-function:

$$L(\chi, s) = \sum_{n=1}^{\infty} \frac{\chi(n)}{n^s},$$

for any complex number s with a real part greater than 1. Using analytic continuation, this function can be extended to a meromorphic function defined on the entire complex plane.

The generalized Riemann conjecture states that for any Dirichlet character χ and any complex number s with L (χ, s) = 0, it holds that if a real number s is between 0 and 1, then it is, in fact, equal to 1/2.

If the generalized Riemann conjecture is true, then the corollary is fulfilled - an asymptotic estimate for the number of primes located on a geometric progression $kn + l$ on a ray $[2, x)$ is:

$$\pi_l(k,x) = \frac{1}{\varphi(k)} \int_2^x \frac{dt}{\log(t)} + O(x^{1/2} \log(x)). \tag{1.6}$$

Based on the indicated asymptotics of the number of primes located on a geometric progression $kn+l$ on a ray $[2, x)$, a general formula for the asymptotic estimate of the sums of functions of primes on a given geometric progression will be obtained.

Using the asymptotic law of distribution of primes located on a geometric progression, asymptotics of the sums of various functions of primes on a given geometric progression are determined.

Earlier, asymptotic estimates for the sums of functions of prime numbers on a geometric progression were very laborious, not systematized, and were obtained from other considerations [5].

This approach makes it possible to first state the asymptotic law of distribution of prime numbers in a geometric progression and the generalized Riemann conjecture, and then, based on them, derive a general formula for the asymptotic estimate of the sums of functions of prime numbers in a geometric progression. As an example of using the general formula we perform the derivation separate formulas, as done in Chapter 2 of this work.

Necessary and sufficient conditions for existence the asymptotics will be proved in Chapter 3 of the work.

## 2. ASYMPTOTICS OF THE SUMS OF FUNCTIONS OF PRIME NUMBERS ON A GEOMETRIC PROGRESSION

Let $a_m$ is a sequence of real or complex numbers and $f(m)$ is a continuously differentiable function on a ray $[1, x)$, then the Abel summation formula is true [6]:

$$\sum_{p \leq x, p \in km+l} f(p) = \sum_{m=2}^{x} a_m f(m) = A(x) f(x) - \int_2^x A(t) f'(t) dt, \tag{2.1}$$

where $A(x) = \sum_{m=1}^{x} a_m$.

Let's look at the distribution of primes on a geometric progression. Let $a_m = 1$, if the value $km+l$ is a prime number, where $(k,l)=1$, and $a_m = 0$ - otherwise. Then $A(x) = \pi_l(k,x)$, where $\pi_l(k,x)$ is the number of primes on a geometric progression $km+l$ on a ray $[2,x)$.

Substituting $A(x) = \pi_l(k,x)$ into formula (2.1) and get:

$$\sum_{p \leq x, p \in km+l} f(p) = \sum_{m=2}^{x} a_m f(m) = \pi_l(k,x) f(x) - \int_2^x \pi_l(k,t) f'(t) dt \quad . \tag{2.2}$$

We will use (2.2) to prove a number of assertions.

First, we use the asymptotic estimate (1.2) to determine the asymptotic behavior of the sums of functions of primes located on a geometric progression $km+l$, where $(k,l)=1$.

Assertion 2.1

If $f(t)$ is a monotone function and has a continuous derivative on the ray $[2,x)$, then for any $k,l$:

$$\sum_{p \leq x, p \in km+l} f(p) = \frac{xf(x)}{\varphi(k)\log(x)} - \frac{1}{\varphi(k)} \int_2^x \frac{tf'(t)dt}{\log(t)} + O\left(\frac{x|f(x)|}{\log^2(x)}\right) + O\left(\int_2^x \frac{t|f'(t)|dt}{\log^2(t)}\right). \tag{2.3}$$

Substituting (1.2) into (2.2) for the proof, and we obtain (2.3).

Let us give examples of using formula (2.3):

1. $\displaystyle\sum_{p \leq x, p \in km+l} 1 = \frac{x}{\varphi(k)\log(x)} + O\left(\frac{x}{\log^2(x)}\right)$, which corresponds to formula (1.2).

2. $\displaystyle\sum_{p \leq x, p \in km+l} \log(p) = \frac{x}{\varphi(k)} + O(x/\log(x))$ \hfill (2.4)

3. $\displaystyle\sum_{p \leq x, p \in km+l} 1/p = \frac{\log\log(x)}{\varphi(k)} + O(1)$ \hfill (2.5)

Now we use a more precise asymptotic estimate (1.1).

Assertion 2.2

If the function $f$ is monotone and has a continuous derivative on the ray $[2, x)$, then for any $k, l$:

$$\sum_{p \leq x, p \in km+l} f(p) = \frac{1}{\varphi(k)} \int_2^x \frac{f(t)dt}{\log(t)} + O\left(\frac{|f(x)|x}{e^{c\log^{1/2}(x)}}\right) + O\left(\int_2^x \frac{t|f'(t)|dt}{e^{c\log^{1/2}(t)}}\right),$$

where $c$ is a constant.

Proof

Substituting (1.1) into (2.2) and we get:

$$\sum_{p \leq x, p \in km+l} f(p) = \frac{f(x)}{\varphi(k)} \int_2^x \frac{dt}{\log(t)} + O\left(\frac{|f(x)|x}{e^{c\log^{1/2}(x)}}\right) - \frac{1}{\varphi(k)} \int_2^x \left(\int_2^t \frac{du}{\log(u)}\right) f'(t)dt + O\left(\int_2^x \frac{t|f'(t)|dt}{e^{c\log^{1/2}(t)}}\right). \quad (2.6)$$

We use the method of integration by parts:

$$\frac{1}{\varphi(k)} \int_2^x \left(\int_2^t \frac{du}{\log(u)}\right) f'(t)dt = \frac{f(x)}{\varphi(k)} \int_2^x \frac{du}{\log(u)} - \frac{1}{\varphi(k)} \int_2^x \frac{f(t)dt}{\log(t)} \quad (2.7)$$

Substituting (2.7) into (2.6) and get:

$$\sum_{p \leq x, p \in km+l} f(p) = \frac{1}{\varphi(k)} \int_2^x \frac{f(t)dt}{\log(t)} + O\left(\frac{|f(x)|x}{e^{c\log^{1/2}(x)}}\right) + O\left(\int_2^x \frac{t|f'(t)|dt}{e^{c\log^{1/2}(t)}}\right), \quad (2.8)$$

which corresponds to assertion 2.2.

Let's take a look at examples of use (2.8):

1. $\sum_{p \leq x, p \in km+l} 1 = \frac{1}{\varphi(k)} \int_2^x \frac{dt}{\log(t)} + O\left(\frac{x}{e^{c\log^{1/2}(x)}}\right)$, which corresponds to (1.1).

2. $\sum_{p \leq x, p \in km+l} \log(p) = \frac{1}{\varphi(k)} \int_2^x \frac{\log(t)dt}{\log(t)} + O\left(\frac{\log(x)x}{e^{c\log^{1/2}(x)}}\right) + O\left(\int_2^x \frac{dt}{e^{c\log^{1/2}(t)}}\right) = \frac{x}{\varphi(k)} + O\left(\frac{x\log(x)}{e^{c\log^{1/2}(x)}}\right).$ (2.9)

Compare (2.9) and (2.4).

3. $\sum_{p \leq x, p \in km+l} \frac{\log(p)}{p} = \frac{1}{\varphi(k)} \int_2^x \frac{\log(t)dt}{t\log(t)} + O\left(\frac{\log(x)x}{xe^{c\log^{1/2}(x)}}\right) + O\left(\int_2^x \frac{tdt}{t^2 e^{c\log^{1/2}(t)}}\right) + O\left(\int_2^x \frac{t\log(t)dt}{t^2 e^{c\log^{1/2}(t)}}\right) = \frac{\log(x)}{\varphi(k)} + O(1)$

4. $$\sum_{p \leq x, p \in km+l} p^{\alpha} = \frac{1}{\varphi(k)} \int_2^x \frac{t^{\alpha} dt}{\log(t)} + O(\frac{x^{\alpha+1}}{e^{c\log^{1/2}(x)}}) + O(\int_2^x \frac{t^{\alpha} dt}{e^{c\log^{1/2}(t)}}) = \frac{1}{\varphi(k)} \int_2^x \frac{t^{\alpha} dt}{\log(t)} + O(\frac{x^{\alpha+1}}{e^{c\log^{1/2}(x)}}),$$

where $\alpha > -1$.

Assertion 2.3

If the function $f$ is monotone and has a continuous derivative on the ray $[2, x)$, then if the generalized Riemann conjecture is true, for any $k, l$ is fulfilled:

$$\sum_{p \leq x, p \in km+l} f(p) = \frac{1}{\varphi(k)} \int_2^x \frac{f(t)dt}{\log(t)} + O(|f(x)| x^{1/2} \log(x)) + O(\int_2^x |f'(t)| t^{1/2} \log(t) dt). \qquad (2.10)$$

The proof is similar to assertion 2.2, only the value (1.6) is substituted in (1.2).

The asymptotic estimate for the remainder term (2.10) is the best. Somebody can be verified by comparing the remainder terms in assertion 2.3 with assertion 2.1 and 2.2.

Let us show that the leading terms of the asymptotics (2.3) and (2.8) coincide, and the differences are in the remainders.

$$\sum_{p \leq x, p \in km+l} f(p) = \frac{1}{\varphi(k)} \int_2^x \frac{f(t)dt}{\log(t)} + O(\frac{|f(x)| x}{e^{c\log^{1/2}(x)}}) + O(\int_2^x \frac{t|f'(t)|dt}{e^{c\log^{1/2}(t)}})$$

We transform the term in formula (2.8):

$$\int_2^x \frac{f(t)dt}{\log(t)} = \frac{f(x)x}{\log(x)} - \frac{2f(2)}{\log(2)} - \int_2^x \frac{tf'(t)dt}{\log(t)} + \int_2^x \frac{f(t)dt}{\log(t)}. \qquad (2.11)$$

Substitute (2.11) in (2.8):

$$\sum_{p \leq x, p \in km+l} f(p) = \frac{f(x)x}{\varphi(k)\log(x)} + O(1) - \frac{1}{\varphi(k)} \int_2^n \frac{tf'(t)dt}{\log(t)} + \frac{1}{\varphi(k)} \int_2^n \frac{f(t)dt}{\log(t)} + O(\frac{|f(x)x|}{e^{c\log^{1/2}(x)}}) + O(\int_2^x \frac{t|f'(t)|dt}{e^{c\log^{1/2}(t)}}). (2.12)$$

Let us compare expression (2.12) with the expression for asymptotic (2.3). The first three terms are the same. The discrepancy is only in the residual terms. The remainder terms of expression (2.12) give a more accurate estimate.

## 3. NECESSARY AND SUFFICIENT CONDITIONS FOR THE EXISTENCE OF THESE ASYMPTOTICS

Let $a_m = 1$, if $km + l$ is a prime number, where $(k, l) = 1$, and $a_m = 0$ otherwise.

Let $b_1 = 0, b_m = \dfrac{1}{\varphi(k) \ln m}$, where $m$ is a natural number.

Let's denote: $A(n) = \sum\limits_{m \leq n} a_m = \sum\limits_{m \leq n} a_m f(m)$ and $B(n) = \sum\limits_{m \leq n} b_m = \sum\limits_{m \leq n} b_m f(m)$.

It is required that:

$$\lim_{n \to \infty} \frac{A(n)}{B(n)} = 1. \tag{3.1}$$

We write (3.1) in the form:

$$\lim_{n \to \infty} \frac{\sum\limits_{k=1}^{n} a_k f(k)}{\sum\limits_{k=1}^{n} b_k f(k)} = 1. \tag{3.2}$$

Assertion 3.1

When the conditions are met:

1. $\lim\limits_{n \to \infty} \dfrac{\int_1^n B(t) f'(t) dt}{B(n) f(n)}$ is not equal to 1.

2. $f(x)$ - monotonous and $f'(x) \neq 0$.

3. $\lim\limits_{n \to \infty} \int_1^n B(t) f'(t) dt = \pm \infty$.

It is performed:

$$\lim_{n \to \infty} \frac{\sum\limits_{k=1}^{n} a_k f(k)}{\sum\limits_{k=1}^{n} b_k f(k)} = 1$$

Proof

Having in mind the Abel summation formula:

$$\sum_{m=1}^{n} a_m f(m) = A(n)f(n) - \int_1^n A(t)f'(t)dt, \quad \sum_{m=1}^{n} b_m f(m) = B(n)f(n) - \int_1^n B(t)f'(t)dt. \quad (3.3)$$

Substituting expressions (3.3) into (3.2), we obtain:

$$\frac{\sum_{k=1}^{n} a_k f(k)}{\sum_{k=1}^{n} b_k f(k)} = \frac{A(n)f(n) - \int_1^n A(t)f'(t)dt}{B(n)f(n) - \int_1^n B(t)f'(t)dt} = \frac{A(n)}{B(n)} \cdot \frac{1 - \frac{\int_1^n A(t)f'(t)dt}{A(n)f(n)}}{1 - \frac{\int_1^n B(t)f'(t)dt}{B(n)f(n)}}. \quad (3.4)$$

Let the above sufficient conditions are true:

1. $\lim_{n\to\infty} \frac{\int_1^n B(t)f'(t)dt}{B(n)f(n)}$ is not equal to 1.

2. $f(x)$ - monotonous and $f'(x) \neq 0$.

3. $\lim_{n\to\infty} \int_1^n B(t)f'(t)dt = \pm\infty$.

Let us show that under the indicated conditions:

$$\lim_{n\to\infty} \frac{\int_1^n A(t)f'(t)dt}{A(n)f(n)} = \lim_{n\to\infty} \frac{\int_1^n B(t)f'(t)dt}{B(n)f(n)}. \quad (3.5)$$

Since, using the L'Hôpital rule, we get:

$$\lim_{n\to\infty} \frac{\frac{\int_1^n A(t)f'(t)dt}{A(n)f(n)}}{\frac{\int_1^n B(t)f'(t)dt}{B(n)f(n)}} = \lim_{n\to\infty} \frac{B(n)}{A(n)} \cdot \lim_{n\to\infty} \frac{\int_1^n A(t)f'(t)dt}{\int_1^n B(t)f'(t)dt} = \lim_{n\to\infty} \frac{A(n)f'(n)}{B(n)f'(n)} = 1, \quad (3.6)$$

which corresponds to (3.5).

Then, based on (3.1) and (3.6), we obtain:

$$\lim_{n\to\infty}\frac{\sum_{k=1}^{n}a_k f(k)}{\sum_{k=1}^{n}b_k f(k)}=\lim_{n\to\infty}\frac{A(n)}{B(n)}\cdot\frac{1-\dfrac{\int_1^n A(t)f'(t)dt}{A(n)f(n)}}{1-\dfrac{\int_1^n B(t)f'(t)dt}{B(n)f(n)}}=1,$$

which corresponds to the assertion.

Corollary 3.2

Conditions (1), (3) in Assertion 3.1 correspond to:

1. $\displaystyle\lim_{n\to\infty}\frac{\int_2^n \dfrac{tf'(t)}{\log(t)}dt}{\dfrac{nf(n)}{\log(n)}}$ is not equal to 1.

3. $\displaystyle\lim_{n\to\infty}\int_2^n \frac{tf'(t)}{\log(t)}dt = \pm\infty$.

Proof

Let's find the value:

$$B(n)=\sum_{m\leq n}b_m = \frac{1}{\varphi(k)}\sum_{m\leq n}\frac{1}{\ln m} = \frac{n}{\varphi(k)\ln(n)}(1+o(1)). \tag{3.7}$$

Substituting (3.7) into conditions 3 of Assertion 3.1 and we obtain:

$$\lim_{n\to\infty}\int_1^n B(t)f'(t)dt = \frac{1}{\varphi(k)}\int_2^n \frac{tf'(t)dt}{\log t}(1+o(1)). \tag{3.8}$$

Having in mind (3.8) and if $\displaystyle\lim_{n\to\infty}\int_2^n \frac{tf'(t)}{\log(t)}dt=\pm\infty$ the condition (3) is true.

Now we substitute (3.7) into condition 1 of Assertion 3.1 and obtain:

$$\lim_{n\to\infty}\frac{\int_1^n B(t)f'(t)dt}{B(n)f(n)} = \lim_{n\to\infty}\frac{\dfrac{1}{\varphi(k)}\int_2^n\dfrac{tf'(t)dt}{\log t}(1+o(1))}{\dfrac{nf(n)}{\varphi(k)\log(n)}(1+o(1))} = \lim_{n\to\infty}\frac{\int_2^n\dfrac{tf'(t)dt}{\log t}}{\dfrac{nf(n)}{\log(n)}} \quad\text{not equal to 1.}$$

The conditions of Assertion 3.1 and Corollary 3.2 are sufficient for the indicated asymptotics to hold.

Let's consider a special case of Assertion 3.1 and Corollary 3.2, when the function $f$ is monotonically increasing and tends to infinity.

Assertion 3.3

Let $f(n)$ - increase, i.e. $f'(n) > 0$, $\lim_{n \to \infty} f(n) = \infty$, $f'(n)$ is a continuous function and $\lim_{n \to \infty} \frac{f(n)}{\log(n) f'(n)} \neq 0$. Then all the conditions of Assertion 3.1 and Corollary 3.2 are true.

Proof

Let us check the first condition of Corollary 3.2 using L'Hôpital:

$$\lim_{n \to \infty} \frac{\int_2^n \frac{tf'(t)}{\log(t)} dt}{\frac{nf(n)}{\log(n)}} = \lim_{n \to \infty} \frac{\frac{nf'(n)}{\log(n)}}{(\frac{nf(n)}{\log(n)})'} = \lim_{n \to \infty} \frac{\frac{nf'(n)}{\log(n)}}{\frac{nf'(n)}{\log(n)} + (\frac{n}{\log(n)})' f(n)} = \lim_{n \to \infty} \frac{1}{1 + \lim_{n \to \infty} (\frac{f(n)}{nf'(n)} - \frac{f(n)}{\log(n) f'(n)})} \neq 1. \quad (3.9)$$

It is required $\lim_{n \to \infty} \frac{f(n)}{\log(n) f'(n)} \neq 0$ that (3.9) is fulfilled.

Let us check the second condition of Assertion 3.1:

$f(n)$ is a monotone function and $f'(n) \neq 0$ by the condition of Assertion 3.3.

Let us check the third condition of Corollary 3.2 and show that $\lim_{n \to \infty} \int_2^n \frac{tf'(t)}{\log(t)} dt = \pm \infty$.

$$\int_2^n \frac{tf'(t)}{\log(t)} dt \geq \int_2^n \frac{2f'(t)}{\log(2)} dt = \frac{2}{\log(2)} \int_2^n f'(t) dt = \frac{2}{\log(2)} (f(n) - f(2)) \text{ - increases indefinitely,}$$

like $f(n)$.

Let us show that the conditions of Assertion 3.3 are true for the function $f(n) = n^l, l \geq 0$.

$f(n)$ - increases, $\lim_{n \to \infty} f(n) = \infty$, $f'(n)$ - continuous function,

$$\lim_{n \to \infty} \frac{f(n)}{\log(n) f'(n)} = \lim_{n \to \infty} \frac{n^l}{\log n \, l(n)^{l-1}} = \infty \neq 0.$$

Let us show that the conditions of Assertion 3.3 are not true for the function $f(n) = 2^n$.

$f(n)$ - increases, $\lim_{n \to \infty} f(n) = \infty$, $f'(n)$ - continuous function,

$$\lim_{n \to \infty} \frac{f(n)}{\log(n) f'(n)} = \lim_{n \to \infty} \frac{2^n}{\log n \, 2^n \log 2} = 0.$$

Let us now consider a necessary condition for the fulfillment of the indicated asymptotics.

We will use the notation of Assertion 3.1.

Assertion 3.4

Let:

$$\lim_{n \to \infty} \frac{\sum_{m=1}^{n} a_m f(m)}{\sum_{m=1}^{n} b_m f(m)} = 1$$

Then, for $p \to \infty$ ($p$ is a prime number):

$$\left| \frac{f(p)}{\sum_{m=1}^{p} b_m f(m)} \right| \to 0. \tag{3.10}$$

Proof

If:

$$\lim_{n \to \infty} \frac{\sum_{m=1}^{n} a_m f(m)}{\sum_{m=1}^{n} b_m f(m)} = 1,$$

then:

$$\lim_{n \to \infty} \left| \frac{\sum_{m=1}^{n} a_m f(m)}{\sum_{m=1}^{n} b_m f(m)} - \frac{\sum_{m=1}^{n-1} a_m f(m)}{\sum_{m=1}^{n-1} b_m f(m)} \right| = 0, \tag{3.11}$$

Let's take $n = p$ where $p$ is a prime number.

Then, based on (3.11):

$$\lim_{n \to \infty} \left| \frac{\sum_{m=1}^{p} a_m f(m)}{\sum_{m=1}^{p} b_m f(m)} - \frac{\sum_{m=1}^{p-1} a_m f(m)}{\sum_{m=1}^{p-1} b_m f(m)} \right| = 0 \qquad (3.12)$$

Let us denote $m = p$ by value $a_m f(m) = a_p f(p)$.

We use that $a_p = 1$ and get:

$$\sum_{m=1}^{p} a_m f(m) = a_p f(p) + \sum_{m=1}^{p-1} a_m f(m) = f(p) + \sum_{m=1}^{p-1} a_m f(m). \qquad (3.13)$$

We transform (3.12) and, taking into account (3.13) and obtain:

$$\left| \frac{\sum_{k=1}^{p} a_m f(m)}{\sum_{k=1}^{p} b_m f(m)} - \frac{\sum_{k=1}^{p-1} a_m f(m)}{\sum_{k=1}^{p-1} b_m f(m)} \right| = \left| \frac{a_p f(p) \sum_{k=1}^{p-1} b_m f(m) - b_p f(p) \sum_{k=1}^{p-1} a_m f(m)}{\sum_{m=1}^{p} b_m f(m) \sum_{k=1}^{p-1} b_m f(m)} \right| = |f(p)| \left| \frac{1 - b_p \frac{\sum_{m=1}^{p-1} a_m f(m)}{\sum_{m=1}^{p-1} b_m f(m)}}{\sum_{m=1}^{p} b_m f(m)} \right|$$

Since for $p \to \infty$, value $b_p \to 0$ and $\lim_{n \to \infty} \frac{\sum_{m=1}^{n} a_m f(m)}{\sum_{m=1}^{n} b_m f(m)} = 1$, then the following is performed:

$$\left| \frac{f(p)}{\sum_{m=1}^{p} b_m f(m)} \right| \to 0,$$

which corresponds to (3.10).

Let's consider examples of fulfilling the necessary condition:

1. $f(p) = \log p$. Since $\sum_{m=1}^{p} b_m f(p) = \frac{1}{\varphi(k)} \sum_{m=2}^{p} \frac{\log m}{\log m} = \frac{p}{2\varphi(k)}$, then

$$\left|\frac{f(p)}{\sum_{m=1}^{p} b_m f(p)}\right| = \left|\frac{\log p}{\frac{p}{2\varphi(k)}}\right| \mapsto 0 \text{ for } p \to \infty.$$

2. $f(p) = p^l, l > 0$. Since $\sum_{m=1}^{p} b_m f(p) = \frac{1}{\varphi(k)} \sum_{m=2}^{p} \frac{m^l}{\log m} = \frac{p^{l+1}}{\varphi(k) \log(p)} (1 + o(1))$.

Then $\left|\dfrac{f(p)}{\sum_{m=1}^{p} b_m f(p)}\right| = \left|\dfrac{p^l}{\dfrac{p^{l+1}}{\varphi(k)\log(p)}(1+o(1))}\right| \mapsto 0$ for $p \to \infty$.

Let's look at an example of a prerequisite failure:

3. $f(p) = \dfrac{1}{p^2}$. Since $\sum_{m=1}^{p} b_m f(p) = \dfrac{1}{\varphi(k)} \sum_{m=2}^{p} \dfrac{1}{m^2 \log m} = \dfrac{C}{\varphi(k)}$, where C is constant.

Then $\lim_{p\to\infty} \left|\dfrac{f(p)}{\sum_{m=1}^{p} b_m f(p)}\right| = \lim_{p\to\infty} \left|\dfrac{\frac{1}{p^2}}{\frac{C}{\varphi(k)}}\right| \neq 0$.

4. CONCLUSION AND SUGGESTIONS FOR FURTHER WORK

The next article will continue to study the behavior of some sums of functions of primes.

5. ACKNOWLEDGEMENTS

Thanks to everyone who has contributed to the discussion of this paper. I am grateful to everyone who expressed their suggestions and comments in the course of this work.

References

1. Volfson V.L. Asymptotics of the sums of functions of prime numbers, Applied Physics and Mathematics No. 4, 2020, pp. 23-29.

2. Charles de la Vallée Poussin. Recherces analytiques sur la théorie des nombres premiers. Ann. Soc. Sci. Bruxells, 1897.

3. IM Vinogradov, A new estimate of a function, Izv., Academy of Sciences of the USSR, ser., Mat., No. 2, t 22 (1958), 161-164.

4. E.K. Tichmarsh, The Theory of the Riemann Zeta Function, 1953, Trans. from English, - M.: IL -409 pp.

5. K. Prakhar, Distribution of prime numbers, M, Mir, 1967 -512 p.

6. Apostol, Tom, Introduction to Analytic Number Theory, Undergraduate Texts in Mathematics, Springer-Verlag (1976).